\documentclass[10pt, conference]{IEEEtran}
\IEEEoverridecommandlockouts
\pdfoutput=1 
\makeatletter
 \def\ps@headings{%
 \def\@oddhead{\mbox{}\scriptsize\rightmark \hfil \thepage}%
 \def\@evenhead{\scriptsize\thepage \hfil \leftmark\mbox{}}%
 \def\@oddfoot{}%
 \def\@evenfoot{}}
 \makeatother
\usepackage{amsfonts,amsthm,amssymb,amsmath}
\usepackage{graphicx,multirow,dsfont}
\usepackage[usenames,dvipsnames]{color}

\usepackage{hyperref}

\usepackage{tikz}
\usetikzlibrary{arrows,shapes}

\usepackage[normalem]{ulem}
\newtheorem{defn}{Definition}%[section]
\newtheorem{thm}{Theorem}%[section]
\newtheorem{prop}{Proposition}%[section]
\newcommand{\R}{\mathbb{R}} % real numbers
\newcommand{\pr}{\mathbb{P}} % probability
\newcommand{\E}{\mathbb{E}} % expectation
\newcommand{\dto}{\Rightarrow} % weak convergence
\newcommand{\ds}[1]{\hat{#1}^n} % diffusion scaling
\newcommand{\invar}{\mathcal{I}}
\newcommand{\lyp}{\mathcal{L}}

\usepackage{cite} %Maria's change, based on the IEEE template instructions

\title{Separation of Timescales in a Two-Layered Network}

\author{
\IEEEauthorblockN{Maria Vlasiou\IEEEauthorrefmark{1},
Jiheng Zhang\IEEEauthorrefmark{2},
Bert Zwart\IEEEauthorrefmark{3},
Rob van der Mei\IEEEauthorrefmark{3}}
\IEEEauthorblockA{\IEEEauthorrefmark{1}Department of Mathematics and Computer Science\\
Eindhoven University of Technology, P.O.\ Box 513, 5600 MB, Eindhoven, The Netherlands\\
Email: m.vlasiou@tue.nl}
\IEEEauthorblockA{\IEEEauthorrefmark{2}Department of Industrial Engineering and Logistics Management\\
Hong Kong University of Science and Technology, Hong Kong, S.A.R., China\\
Email: j.zhang@ust.hk}
\IEEEauthorblockA{\IEEEauthorrefmark{3}Centrum Wiskunde \& Informatica\\
Science Park 123, Amsterdam, The Netherlands\\
Email: (bert.zwart, mei)@cwi.nl}
\thanks{The research of Maria Vlasiou and Jiheng Zhang is partly supported by two grants from the `Joint Research Scheme' program, sponsored by the Netherlands Organization of Scientific Research (NWO) and the Research Grants Council of Hong Kong (RGC) through projects 649.000.005 and D-HK007/11T, respectively. The research of Bert Zwart is partly supported by an NWO VIDI grant and an IBM faculty award.}
}

\usepackage{ifpdf}
\ifpdf
\else
\renewcommand{\includegraphics}[1][]{\url}
\fi
\begin{document}

\maketitle

\begin{abstract}
We investigate a computer network consisting of two layers occurring in, for example, application servers. The first layer incorporates the arrival of jobs at a network of multi-server nodes, which we model as a many-server Jackson network.  At the second layer, active servers at these nodes act now as customers who are served by a common CPU. Our main result shows a separation of time scales in heavy traffic: the main source of randomness occurs at the (aggregate) CPU layer; the interactions between different types of nodes at the other layer is shown to converge to a fixed point at a faster time scale; this also yields a state-space collapse property. Apart from these fundamental insights, we also obtain an explicit approximation for the joint law of the number of jobs in the system, which is provably accurate for heavily loaded systems and performs numerically well for moderately loaded systems. The obtained results for the model under consideration can be applied to thread-pool dimensioning in application servers, while the technique seems applicable to other layered systems too.
%Keywords: layered queuing networks, web server modeling, diffusion approximation, state-space collapse, heavy traffic.
\end{abstract}

\section{Introduction}

Communication networks need to support a growing diversity and heterogeneity in applications. Examples are web-based multi-tiered system architectures, with a client tier to provide an interface to end users, a business logic tier to coordinate information retrieval and processing, and a data tier with legacy systems to store and access customer data. In such environments, different applications compete for access to shared infrastructure resources, both at the software level (e.g., mutex and database locks, thread-pools) and at the hardware level (e.g., bandwidth, processing power, disk access). Thus, the performance of such applications is determined by the interplay of software and hardware contention. For background, see \cite{WBM2009,Mei01, Cardellini}.

In particular, in situations where web pages are created on-the-fly (think of making a reservation online), the benefits of caching are limited and sizes of web pages are unknown, and there is usually ample core network bandwidth available at reasonable prices. Consequently, the bottleneck in user-level performance can shift from the network interface to the application server, and implementing size-based scheduling policies becomes hard, contrary to the situation considered in \cite{Crovella,Harchol}.

Application servers usually implement a number of thread-pools; a thread is software that can perform a specific type of sub-transaction. Consider for example the web-server performance model proposed in \cite{Mei01}. Each HTTP request that requires server-side scripting (e.g., CGI or ASP scripts, or Java servlets) consists of two subsequent phases: a document-retrieval phase, and a script processing phase. To this end, the web server implements two thread-pools, one performing the first phase of processing, and the other performing the second phase of processing. The model consists of a tandem of two multi-server queues, where servers at queue 1 represent the phase-1 threads, and the servers at queue 2 represent phase-2 threads. A particular feature of this model is that at all times the active threads share a common Central Processing Unit (CPU) in a Processor-Sharing (PS) fashion; cf.\ \cite{Jonckheere10,Wei04}. Alternatively, one can think of scheduling jobs in data centers, where different parts of a job are taken care of by a different thread-pool.

Motivated by this, we study a relatively simple, but nontrivial two-layered network. An informal model description is as follows. The first layer models the processing of jobs by a network of nodes. Each node consists of several servers and, therefore,  it looks like a (generalized) Jackson network consisting of many-server queues. The servers in this network act as customers in a second layer, in the sense that they are served by a single CPU in a PS fashion. A detailed model description is provided in Section~\ref{sec:model-description}.

Variations of the above model have been investigated in several papers in the literature, but apart from stability analysis \cite{Jonckheere10}, a rigorous analysis of this layered network has been lacking. The same can be said about other literature on layered networks. Only a limited number of papers focus on the performance of multi-layered queuing networks. A fundamental paper is Rolia and Sevcik \cite{Rolia95}, who propose the Method of Layers, i.e., a closed queuing-network based model for the responsiveness of client-server applications, explicitly taking into account both software and hardware contention. Another fundamental contribution is presented by Woodside et al.\ \cite{Woodside95}, who propose the so-called Stochastic Rendezvous Network model to analyze the performance of application software with client-server synchronization. The contributions presented in \cite{Rolia95} and \cite{Woodside95} are often referred to as Layered Queuing Models. A common drawback of multi-layered queuing models is that exact analysis is primarily restricted to special cases, and numerical algorithms are typically required to obtain performance measures of interest (see for example \cite{Woodside95}). Although such methods are important, it is also valuable to look at layered systems from a more qualitative point of view, which we do in this paper by considering the system under critical load.

The most simple example of the layered systems we consider is the case where the first layer consists of a single node. In this case, the model reduces to the so-called limited processor sharing (LPS) queue. Recently, there has been considerable interest in the analysis of LPS systems. Avi-Itzhak and Halfin \cite{Avi-ItzhakHalfin1988} propose an approximation for the mean response time.  A computational analysis based on matrix geometric methods is performed in Zhang and Lipsky \cite{ZhaLi06a,ZhaLi07b}.  Some stochastic ordering results are derived in Nuyens and van der Weij \cite{NuyensWeij07}. Large deviation results are presented in Nair {\em et al.} \cite{Nair10}, and these results are also applied to show that LPS provides robust performance across a range of both heavy-tailed and light-tailed job sizes, as it combines the attractive properties of a guaranteed service rate of FIFO and the possibility of overtaking offered by PS.

The work on LPS that is most relevant for this study is the work of Zhang, Dai and Zwart \cite{ZhangZwart2008,ZDZ2009,ZDZ2011} who study the stochastic processes that underlie the LPS queue in the heavy-traffic regime, i.e.\ an asymptotic regime where the traffic intensity converges to 1.  The setting is rather general, allowing the inter-arrival and service times to have general distributions. Fluid and diffusion limits are derived, leading to a heavy-traffic analysis of the steady-state distribution of LPS, showing that the approximation by Avi-Itzhak and Halfin \cite{Avi-ItzhakHalfin1988} is asymptotically accurate in heavy traffic.

In the present paper, we perform an analysis similar to the one performed in \cite{ZhangZwart2008,ZDZ2009,ZDZ2011}. Under the assumption that job sizes are exponentially distributed, we expand the work in \cite{ZhangZwart2008,ZDZ2009,ZDZ2011} from the single node case to networks. Moreover, based on our mathematical results, we propose an extension to general job sizes.

We analyze the system as it approaches heavy traffic. Under the assumption that there is a single bottleneck (an exact definition of bottleneck is given later), we derive explicit results for the joint distribution of the number of jobs in the system by proving a diffusion limit theorem. This limit theorem does not only yield explicit approximations but yields also useful insights: if we look at the system from the CPU layer, we can aggregate the whole system since the total workload acts as if we were dealing with a single server queue. However, information regarding the interaction of several types of customers at the other layer would then be lost. It turns out, nonetheless, that those interactions take place at a much faster time scale in heavy traffic, and that the number of users of all types converge instantaneously to a piece-wise linear function of the number of users at the bottleneck. This separation of time scales property is shown to imply that in heavy traffic, the joint queue length vector can be written as a deterministic function of the total workload as seen from the CPU layer. Such a property is known as \textsl{state-space collapse} (SSC) in the stochastic network literature.

Thus, our methodological contribution is that it is possible to rigorously establish a separation of time scales property in heavy traffic in an important class of layered networks, which makes these layered networks tractable. Although we focus on the Markovian case, we believe that such properties hold more generally as well; we provide some physical and numerical arguments to support this claim. The result on separation of time scales result essentially implies that the main source of randomness in heavy traffic can be observed at the CPU layer, thus making performance analysis much more tractable. Apart from supporting these claims by theorems, some numerical experiments suggest that the resulting approximations perform well. The results in our paper may be useful to create design rules, for example to dimension thread-pools. Some first efforts using heuristic approximations were proposed in \cite{Wei04}.

The paper is organized as follows. We provide a detailed model description in Section~\ref{sec:model-description}. In Section~\ref{sec:fluid-analys-invar} we propose a fluid model for our two-layered system. We use this fluid model to analyze how users of different types interact if the system is in heavy traffic. In doing so, we construct a Lyapounov function which we use to show that the user population converges uniformly to a fixed point that is uniquely defined through the total workload. The fluid model also helps understand which stations will be bottlenecks. Section~\ref{sec:state-space-collapse} contains our main results, namely a process limit theorem for the customer population process. A heavy-traffic approximation of the steady-state distribution of the customer population is proposed in Section~\ref{s:Steady-state performance approximations}. Section~\ref{sec:extens-gener-job} presents an extension to general service times based on physical arguments, and some numerical validation by comparing the proposed approximations with simulation results. Concluding remarks can be found in Section~\ref{s:Concluding remarks}.

%\textcolor{Mulberry}{This part of the intro needs to be fixed. Promised sections that won't be written, and simulation %ending in?}

%Currently, in-depth understanding of the behavior
%of multi-layered performance models is lacking. Motivated by this, we study a simple
%but non-trivial multi-layered queuing model, a two-layered tandem of two multi-server queues,
%where each of the active servers share an underlying resource in a processor-sharing (PS) fashion.

\section{Model description}
\label{sec:model-description}

The purpose of this section is to give a formal model description. We adopt the convention that all vectors are column vectors, and use $a^T$ to denote the transpose of a vector or matrix. For two vectors $x$, $y$ we denote $xy$ to be the vector consisting of elements $x_iy_i$.  Furthermore, $I$ is the identity matrix, $e$ is the vector consisting of $1$'s, and $e_i$ is the vector whose $i$th element is 1 and the rest are all 0. Last, $(x)^+ = \max\{0,x\}$ and $a \wedge b = \min\{a,b\}$.

We consider a network with $J$ nodes. Jobs arrive at node $i\in\{1,\ldots,J\}$ according to a Poisson process with rate $\lambda_i$. Jobs have a random amount of service requirement, which is exponentially distributed with rate $\mu_i$.  Node~$i$ has $K_i$ servers, which allows for parallel processing for the first $K_i$ jobs at the node. Customers move between queues according to a substochastic routing matrix $P$ of dimension $J$.  As in the case of regular queuing networks, we need to introduce the total arrival rates of jobs to station $i$ (i.e.\ including the external arrival rate $\lambda_i$ and internal arrivals from other nodes), which are denoted by $\gamma_i$.  The arrival rates $\gamma_i$ can be found as the unique solution to a system of linear traffic equations. Let $\gamma $ be the vector of elements $\gamma_i$, i.e.\ $\gamma =(\gamma_i)$, and similarly $\lambda =(\lambda _i)$. Then, in vector form, the traffic equation(s) can be written as
\begin{equation*}
  \label{eq:traffic}
  \gamma = \lambda + P^T\gamma.
\end{equation*}
Throughout the paper, we need to assume that $I-P^T$ is invertible, as is usual for open Jackson networks, which leads to the unique solution $\gamma=(I-P^T)^{-1}\lambda$. All active servers interact since they share a CPU working at rate 1. In other words, from the viewpoint of the individual nodes, we have a multi-dimensional Markovian queuing network where jobs of type $i$ are served at rate
\begin{equation}
  \label{eq:R-def}
  R_i(x) := \frac{\min \{x_i,K_i\}}{\sum_j\min\{x_j,K_j\}},
\end{equation}
with $x_i, i=1,\ldots,J$ being the number of customers of type $i$ that are currently in the system. This is consistent with the fact that $\min\{x_j,K_j\}$ is the number of busy servers of type $j$, and all busy servers share the common CPU according to the PS discipline.

It can be useful to view the system from the CPU layer (i.e.\ the second layer), since there is a connection with an $M/PH/1$ queue which we now describe: users arrive at rate $\lambda^o= \sum_i\lambda_i$ and start their service at node $i$ with probability $a_i=\lambda_i/\lambda^o$. Define $a_0=0$, $p_{00}=1$, and for $i\geq 1$, $p_{0i}=0$ and $p_{i0} = 1-\sum_{j}p_{ij}$. Observe that the total service requirement of a job is the time to absorption in state $0$ of a continuous-time Markov chain with initial distribution $(a_i)$, where the time in state $i$ is exponentially distributed with rate $\mu_i$, after which one jumps to state $j$ with probability $p_{ij}$. Thus, the total service requirement $S$ of an arbitrary customer has a phase-type distribution with parameters $(a, \mu, P)$, with $\mu=(\mu_i)$. We also denote by $\beta_i=1/\mu_i$ and $\beta_i^{(2)}= 2/\mu_i^2$ the first and second moment of service requirements at node $i$. The corresponding vectors are denoted by $\beta$ and $\beta^{(2)}$.

It is possible to compute the first two moments of this distribution by using standard methods (see for e.g.\ \cite{Asmussen2003} and references therein).  Let $T_i$ be the total service requirement of each user waiting to be served at node~$i$. This includes their immediate service at node~$i$ and all the future services  due to routing. Denote by $\tau_i$ and $\tau_i^{(2)}$ the first and second moment of $T_i$, and let $\tau,\tau^{2}$ be the corresponding vectors. Then $ \tau = (I-P)^{-1}\beta$  and \begin{equation*}
  \tau_i^{(2)} = \beta_i^{(2)}+ \sum_{j} p_{ij} (2\beta_i \tau_j+\tau_j^{(2)}).
\end{equation*}
In vector notation, this becomes
\begin{equation*}
    \tau^{(2)} = (I-P)^{-1}\left(\beta^{(2)} + 2 \beta (P \tau)\right).
\end{equation*}
Notice that the expressions for $\tau$ and $\tau^{(2)}$ are still valid if the service requirement of a user at node $i$ is not exponential but generally distributed. In that case, the total service requirement is simply the time to absorption of a semi-Markov process. We need this interpretation in Section~\ref{sec:extens-gener-job}. Of course, in that case, it no longer holds that $\beta_i^{(2)}= 2/\mu_i^2$.

We can compute the first and the second moment of the total service requirement $S$, obtaining
$$
E[S] = a^T \tau \mbox{ and } E[S^2] = a^T \tau^{(2)}.
$$
It is also clear from the $M/PH/1$ interpretation that the global stability condition of the system is $E[S]\sum_{i}\lambda_i<1$, or equivalently

\begin{equation*}
%  \label{eq:rho-def}
  \rho := \lambda^T (I-P)^{-1} \beta = \beta^T \gamma<1.
\end{equation*}
We also define $\rho_i = \beta_i\gamma_i=\gamma_i/\mu_i$. Observe that $\rho = \sum_i \rho_i$.

\noindent
{\bf Example:}
We are particularly interested in the simple two-node tandem case ($J=2$), where all users first enter station 1 ($\lambda_2=0$), then move from station 1 to station 2 ($p_{12}=1$) and then leave ($p_{20}=1$). In this case $\gamma_1=\gamma_2=\lambda_1$, $E[S] = 1/\mu_1+1/\mu_2$, and
\[
E[S^2] = 2/\mu_1^2 + 2/(\mu_1\mu_2)+ 2/\mu_2^2.
\]

\begin{figure}[h]\label{fig:LPS-Tandem}
  \begin{center}
    \includegraphics[width=0.48\textwidth]{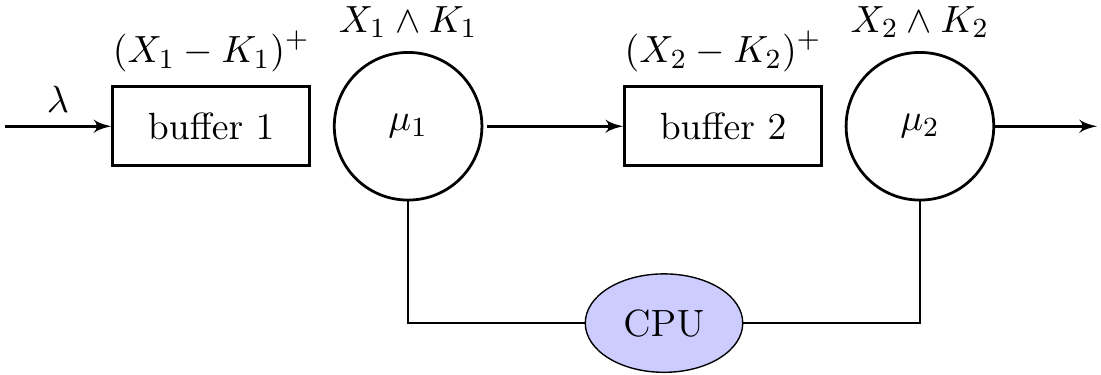}
  \end{center}
  \caption{LPS queues in tandem with a shared CPU.}
\end{figure}
We now investigate the system under critical load, i.e. when $\rho$ is (close to) 1. To this end, we first develop and analyze a critical fluid model in the next section.

\section{Fluid analysis and invariant points}
\label{sec:fluid-analys-invar}

In this section we propose a fluid model for our layered system under the assumption of critical loading, i.e.\ $\rho=1$, or equivalently,
\begin{equation}
  \label{eq:cond-balance}
  \sum_i \frac{\gamma_i}{\mu_i}=1.
\end{equation}
In the sequel, we establish that in this scenario the workload will stay constant, and that the queue length vector will converge to an invariant point. We also characterize the set of invariant points, and show this set is one-dimensional under the assumption that there is a unique bottleneck.

Our fluid model is defined by the following ordinary differential equation (ODE):
\begin{equation}
  \label{eq:X-i}
  \bar X'_i(t)=\lambda_i-\mu_iR_i(\bar X(t))+\sum_{j=1}^Jp_{j,i}\mu_jR_j(\bar X(t)).
\end{equation}
Here, $R_j(x)$ is defined in the same way as in the original stochastic model, cf.\ \eqref{eq:R-def}. Moreover, $\bar{X}_i(t)$ can be interpreted as a fluid approximation of the number of jobs at time $t$, after an appropriate normalization of time and space. We avoid a technical discussion on fluid approximations and refer to \cite{MMR1998} for background.

In particular, for our system it is possible to show the following. Consider a sequence of `virtual' systems indexed by $n$, where the number of servers at node $i$ is equal to $nK_i$ (we call such systems virtual since there exists only one, rather than a sequence of real systems). The network structure, represented by the routing matrix $P$, and the service time at each node is kept fixed. Let $X_i^n(t)$ be the number of type $i$ jobs at time $t$ in the $n$th system. Then it can be shown that $\{X^n_i(nt)/n, t\geq 0, i=1,\ldots,J\}$ converges in the space of functions to $\{\bar X_i(t),t\geq 0, i=1,\ldots,J \}$; see also \cite{MMR1998}. We will not pursue a proof of this fluid limit result here, since it is not our main point. The fluid model we present has a different purpose: it serves as building block for developing a heavy-traffic approximation.

Regarding the scaling of the number of servers $K_i$ to $nK_i$, the skeptical reader should consider that this scaling eventually leads to tractable heavy-traffic approximations in the single-node case as shown in \cite{ZhangZwart2008} and, more importantly, also in the network case as shown later in this paper. In fact, letting the number of servers grow with $n$ is the only way to keep the probability of delay strictly between $0$ and $1$ in heavy traffic. For example, keeping the number of servers fixed would lead to a delay probability of $1$, which is not a very useful approximation for design purposes. In Section~\ref{s:Steady-state performance approximations}, we come back to this limiting procedure, and explain how we can utilize the limit of our sequence of `virtual' systems to obtain performance approximations for the actual system.

Getting back to the fluid model, we can write \eqref{eq:X-i} into vector form
\begin{equation}
  \label{eq:ODE}
  \bar X'=\Psi(\bar X),
\end{equation}
where $\Psi:[0,\infty)^J\to\R^J$ can be represented as
\begin{equation}
  \label{eq:ODE-Phi}
  \Psi(x)=\lambda-\mu R(x)+P^T(\mu R(x)),
\end{equation}
where $R(x)$ is the vector with elements$R_i(x)$ and $\mu R(x)$ indicates a component-wise product, as before.

\begin{thm}[Existence and uniqueness]
For any $\bar X(0)=x\in\R_+^J$, there exist a unique solution to the ODE \eqref{eq:ODE}.
\end{thm}
\begin{proof}
It is clear that each $R_i(x)$ is Lipschitz continuous on $\R_+^J$. So is the linear combination $\Psi(x)$. The result follows from Theorem~VI in Chapter 10 of \cite{Walter1998}.
\end{proof}
Recall that the system is a work-conserving single-server queue when considered at the CPU layer. We now show that this is also the case for our fluid model. We define the workload for the fluid model as follows:
\begin{equation}
  \label{eq:workload-fluid}
  \bar W(t)=\beta^T (1-P^T)^{-1}\bar X(t).
\end{equation}
\begin{prop}\label{prop:workload-conserve}
  For each solution of \eqref{eq:ODE},
$
%  \begin{equation}
%   \label{eq:workload-conserve}
    \bar W(t)=\bar W(0).
%  \end{equation}
$
\end{prop}
\begin{proof}
The proof follows from the computation of the derivative.  From \eqref{eq:ODE},
\begin{align*}
\bar W'(t)&=\beta^T (I-P^T)^{-1}\bar X'(t)\\
    &=\beta^T (I-P^T)^{-1}
    \left(\lambda-\mu R(x)+P^T(\mu R(x))\right)\\
    &=\beta^T \gamma
    - \beta^T (I-P^T)^{-1}(I-P^T) \mu R(\bar X(t))\\
    &=1-\beta^T \mu R(\bar X(t))=1-1=0,
  \end{align*}
where $\beta^T \gamma=1$ is due to critical loading and $\beta^T \mu R(x) =\sum_{i=1}^JR_i(x)=1$ by the definition of $R(x)$ in \eqref{eq:R-def}.
\end{proof}

We now characterize the invariant manifold of the ODE, which is the set of invariant points. A point $x$ is invariant if
\begin{equation}
  \label{eq:balance}
  \mu_iR_i(x)=\gamma_i, \quad i=1,\ldots, J.
\end{equation}
This definition of an invariant point is natural. To see this, observe that the right-hand side of \eqref{eq:balance} represents the total (arrival) rate into node $i$, while the left-hand side of \eqref{eq:balance} can be interpreted as rate out of node $i$, as $R_i(x)$ is the percentage of CPU dedicated to node $i$, thus representing the speed that node $i$ works and $\mu_i^{-1}$ is the service requirement of a job at node $i$.

A crucial notion in the study of invariant points is the notion of bottleneck. It turns out that the following definition is appropriate:
\begin{defn}[Bottleneck]
  \label{def:bottleneck}
  Node $i$ is a bottleneck if $i=\arg\min_j \frac{\mu_jK_j}{\gamma_j}$.
\end{defn}
In this paper, we focus on the case where there is a unique bottleneck. Without loss of generality, we take node~1 as the bottleneck when we investigate the case of a general network; For convenience in numerical experiments and presentation, in the two-node tandem case we may sometimes take node 2 as the bottleneck.

We will now describe the set of invariant points starting from the number of jobs at the bottleneck, i.e.\ $x_1$. There are two cases: if $x_1 < K_1$ then it follows from \eqref{eq:balance} and the definition of $R_i(x)$ that
\[
\mu_ix_i = \gamma_i \sum_j x_j.
\]
Thus, $\sum_jx_j =\mu_1x_1/\gamma_1$, so that $\displaystyle{\mu_ix_i = \gamma_i\tfrac{\mu_1 x_1}{\gamma_1}}$.

In the second case, if $x_1\geq K_1$ then we can write $\displaystyle{\mu_ix_i = \gamma_i\frac{\mu_1 K_1}{\gamma_1}}$. %In general, for $i=2,\ldots, J$, we have that
%\begin{equation*}
%%  \label{eq:inv-one-bottleneck}
% x_i = \frac{\gamma_i}{\mu_i}\frac{\mu_1 (x_1\wedge K_1)}{\gamma_1}.
%\end{equation*}
Thus, the set of invariant points, called the invariant manifold, is
the following:
\begin{equation*}
  \label{eq:inv-single}
  \invar
  =\left\{x\in\R_+^J:
    \frac{\mu_ix_i}{\gamma_i}=\frac{\mu_1(x_1\wedge K_1)}{\gamma_1},
    i=2,\ldots, J
  \right\}.
\end{equation*}
The invariant manifold is illustrated in the following picture for the
two dimensional case.
\begin{figure}[h]
  \centering
  \begin{tikzpicture}[scale=1]
    % Draw axes
    \draw [<->,thick] (0,5) node (yaxis) [above] {$x_1$}
        |- (5,0) node (xaxis) [right] {$x_2$};
    % Draw two intersecting lines
    \draw (0,0) coordinate (a_1) -- (1.6,2) coordinate (a_2);
    \draw (1.6,2) coordinate (b_1) -- (1.6,4.75) coordinate (b_2);
    % \draw[dashed] (0,4.88) -- (2.7111,0);
    % \node at (3.4,0.4) {$s_2(\bar W(t))=K_2$};

    \coordinate (c) at (intersection of a_1--a_2 and b_1--b_2);
    % Draw lines indicating intersection with y and x axis. Here we use
    % the perpendicular coordinate system
    \draw[dashed] (yaxis |- c) node[left] {$K_1$}
        -| (xaxis -| c) node[below] {$\frac{\gamma_2\mu_1K_1}{\mu_2\gamma_1}$};
    \draw[dashed] (2.2, 0.05) -- (2.2,0) node[above] {$K_2$};
    % Draw a dot to indicate intersection point

    % \coordinate (a) at (2.4,0.8);
    % \fill[red] (a) circle (1pt);
    % \draw[blue] (a) -- (1.6,0.8);

    % \coordinate (b) at (1.94,0.9);
    % \fill[red] (b) circle (1pt);
    % \draw[blue] (b) -- (1.44,1.8);

  \end{tikzpicture}
  \caption{Invariant manifold for the 2-dimensional tandem case, where node $1$ is the bottleneck.}
  \label{fig:inv-man-2}
\end{figure}
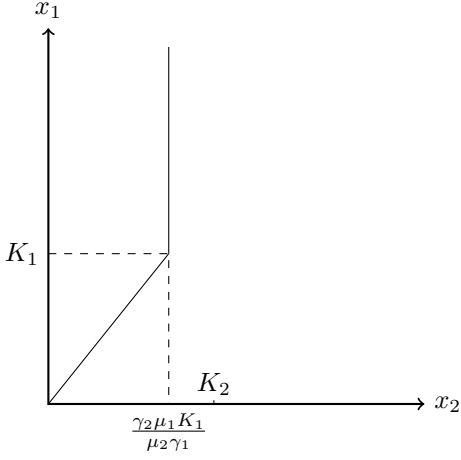

We now conclude by formally showing that our notion of invariant points makes indeed sense.

\begin{prop}
   $\bar X(t)=\bar X(0)$ for all $t\ge 0$ if and only if $\bar X(0)\in\invar$.
\end{prop}

\begin{proof}
The necessity part follows from the above discussion. For sufficiency, it suffices to show that for any $x\in\invar$, $\Psi(x)=0$. Note that by the definition of the invariant manifold, we have that for any $x\in\invar$, $\frac{\mu_ix_i}{\gamma_i}$ is a constant for all $i=2,\ldots, J$. Let $c$ be that constant, i.e.\ $c=\frac{\mu_ix_i}{\gamma_i}$. By \eqref{eq:cond-balance}, $x_1\wedge K_1+\sum_{i=2}^Jx_i=\sum_{i=1}^Jc\frac {\gamma_i}{\mu_i}=c$. Moreover, for $i \geq 2$ and $x\in\invar$, we have that $x_i<K_i$. To see this, observe that $c=\frac{\mu_ix_i}{\gamma_i}\leq \frac{\mu_1x_1}{\gamma_1}$ by the definition of $\invar$, which is less than $\frac{\mu_iK_i}{\gamma_i}$ by the definition of a bottleneck. Thus, we have shown that $R_i(x)=\frac{x_i}{c}$, thus $\mu_iR_i(x)=\gamma_i$. We have
\begin{align*}
    \Psi(x)&= \lambda-\mu R(x)+P^T(\mu R(x))\\
    &=(I-P^T)(\gamma-\mu R(x))=0.\\[-2.2\baselineskip]
\end{align*}
\end{proof}

\subsection{Convergence to invariant points}

We now consider the convergence of the solution $\bar X(t)$ of the ODE \eqref{eq:ODE} to the invariant manifold $\invar$ for any given starting point. Let $x^*$ be the point in the invariant manifold where $x^*_1=K_1$. We easily see that $x^*_i=\beta_i\gamma_i\frac{\mu_1}{\gamma_1}K_1$. Based on this point, we define a critical workload level
\begin{equation*}
  \label{eq:workload-fluid-critical}
  w^*=\beta^T (1-P^T)^{-1}x^*=\beta^T (1-P^T)^{-1} (\beta\gamma) \frac{\mu_1}{\gamma_1}K_1.
\end{equation*}
(Yet another interpretation is  $w^*=\tau^Tx^*$.) This gives rise to a ``critical hyperplane'':
\begin{equation}
  \label{eq:workload-hyperplan}
  \{x: \beta^T (1-P^T)^{-1} x = w^*\}.
\end{equation}

%The idea is simple. The Lyapunov function is constructed based on the distance from the point $x$ to the invariant %manifold depending on whether $x$ is above the critical hyperplane or below.
For any $w\le w^*$, let for $i=1,\ldots, J$
\begin{equation}
  \label{eq:x-dagger-<}
  x^\dagger_i(w) =\frac{\gamma_i w}{\mu_i \beta^T (I-P^T)^{-1} (\beta\gamma) }
  =  \frac{\mu_1K_1}{\gamma_1}\frac{\gamma_i}{\mu_i}\frac{w}{w^*}.
\end{equation}
Note that
\begin{equation*}
R_i(x^\dagger(w))=\frac{x^\dagger_i(w)}{\sum_j x^\dagger_j(w)}
=\frac{\gamma_i/\mu_i}{\sum_j \gamma_j/\mu_j}=\frac{\gamma_i}{\mu_i}.
\end{equation*}
This gives an intuitive explanation that on the invariant manifold, the $R_i$'s, representing outflow of work at station $i$, should be equal to the inflow of work at station $i$. For any $w>w^*$, let
\begin{align}
  \label{eq:x-dagger->-1}
  x^\dagger_1(w)&=K_1
  +\frac{\left(w-w^*\right)}{\beta^T(I-P^T)^{-1}e_1} = K_1+\frac{\left(w-w^*\right)}{\tau_1}\\
  \label{eq:x-dagger->-i}
  x^\dagger_i(w)&=\frac{\mu_1K_1}{\gamma_1}\frac{\gamma_i}{\mu_i},
  \quad i=2,\ldots, J.
\end{align}
It is clear that $x^\dagger(w)$ is the intersection of the workload hyperplane $\bar W(t)=w$ and the invariant manifold. It is now also clear why $w^*$ is called the critical workload level: intuitively, if we are restricted to invariant points, congestion at the bottleneck will occur only if $w>w^*$.

%Depending on whether $w$ is larger than $w^*$ or not, the calculation is different.  We can define a Lyapunov function:
%\begin{equation}
%  \label{eq:Lyapunov}
%  \begin{split}
%    &\lyp(\bar X(t))\\
%    &=\left[
%      \left(\bar X(t)-x^\dagger(w)\right)^T(I-P^T)^{-1}
%      \left(\bar X(t)-x^\dagger(w)\right)
%    \right].
% \end{split}
%\end{equation}
%We now state the following property for the Lyapunov function, which will imply the convergence to the invariant %manifold for the fluid model.
%\begin{prop}
%\label{prop-lyapounov}
%  $\lyp(x)$ is continuous in $x$ and for any $x\notin \invar$, $\lyp(x)>0$; for any
%  $x\in \invar$, $\lyp(x)=0$. When $\bar X(t)\notin\invar$,
%  \begin{equation*}
%    \frac{d}{dt}\lyp(\bar X(t))<0.
%  \end{equation*}
%\end{prop}
%
%A proof of this result can be found in the appendix. The following proposition is an immediate consequence.

The following proposition is the main result of this section. For space limitations, its full proof is postponed to an extended version of this paper.

\begin{prop}[Convergence to the invariant manifold]
  \label{prop:fluid-converge-invar}
  For any solution $\bar X$ to the ODE \eqref{eq:ODE}, we have that
  \begin{equation*}
    \bar X(t)\to x^\dagger(\bar W(0)),\quad \textrm{as }t\to\infty,
  \end{equation*}
  where $x^\dagger$ is as defined by \eqref{eq:x-dagger-<}--\eqref{eq:x-dagger->-i}.
\end{prop}
\begin{proof}
  We can define a Lyapunov function:
  \begin{equation}
    \label{eq:Lyapunov}
    \begin{split}
      &\lyp(\bar X(t))\\
      &=\left[
        \left(\bar X(t)-x^\dagger(w)\right)^T(I-P^T)^{-1}
        \left(\bar X(t)-x^\dagger(w)\right)
      \right].
    \end{split}
  \end{equation}
  It is clear that the function $\lyp(x)$ is continuous in $x$ and for any $x\notin \invar$, $\lyp(x)>0$; for any $x\in \invar$, $\lyp(x)=0$. The result of this proposition will follow immediately if we can show that for any $\bar X(t)\notin\invar$,
  \begin{equation}
    \frac{d}{dt}\lyp(\bar X(t))<0.
  \end{equation}
  Since for any solution $\bar X$ to the ODE \eqref{eq:ODE}, Proposition~\ref{prop:workload-conserve} yields that the workload load $\bar W$ does not change, for $w=\bar W(0)$ we have
  \begin{align}
    \frac{d}{dt}\lyp(\bar X(t))&=
    2\left(\bar X(t)-x^\dagger(w)\right)^T(I-P^T)^{-1}\bar X'(t)\nonumber\\
    &=2\left(\bar X(t)-x^\dagger(w)\right)^T \nonumber\\
    &\quad(I-P^T)^{-1}
      \left[
       \lambda-(I-P^T)\mu R(\bar X(t))
      \right]\nonumber\\
    &=2\left(\bar X(t)-x^\dagger(w)\right)^T\left[
      \gamma-\mu R(\bar X(t))
      \right].\label{eq:lyp-simplified}
  \end{align}
  % It remains to show that for any $i$, if $\bar
  % X_i(t)>x^\dagger_i(w)$, then $\frac{\gamma_i}{\mu_i}\le R_i(\bar
  % X(t))$, and vice versa.

  To simplify the presentation, we focus on the case where the dimension is equal to two, i.e.\ $J=2$.
  Define
  \begin{equation*}
    H_\tau(w)=\left\{x:\tau_1x_1+\tau_2x_2=w\right\}.
  \end{equation*}
  Note that any solution $\bar X$ can only live on $H_\tau(w)$.
  Let $y^\dagger=(x^\dagger_1\wedge K_1,x^\dagger_1\wedge K_1)$.
  It is clear that $R_i(y^\dagger)=\frac{\gamma_i}{\mu_i}$, $i=1,2$.
  If $\bar X_1(t)<x^\dagger_1(w)$, then $\bar X_2(t)>x^\dagger_2(w)$, which follows from  the fact that $\tau_1,\tau_2>0$.  This implies that
%  \begin{align*}
%    \bar X_1(t)\wedge K_1 &\le y^\dagger_1,\\
%    \bar X_2(t)\wedge K_2 &\ge y^\dagger_2.
%  \end{align*}
$$
\bar X_1(t)\wedge K_1 \le y^\dagger_1,\quad \bar X_2(t)\wedge K_2 \ge y^\dagger_2.
$$
Notice that the above two inequalities can not be tight simultaneously. Otherwise, $\bar X(t)$ would be equal to $x^\dagger$. Thus, $R_1(\bar X)\le R_1(x^\dagger(w))$ and $R_2(\bar X)\ge R_2(x^\dagger(w))$ by the definition of $R_i(\cdot)$ in \eqref{eq:R-def}. Again, equality can not hold for both. This implies that $\frac{d}{dt}\lyp(\bar X(t))<0$ according to \eqref{eq:lyp-simplified}. The same argument applies if $\bar X_2(t)<x^\dagger_2(w)$.
\end{proof}

In the following section we show that, as $\rho$ is close to 1, the fluid model is a good approximation of the queue length on a time scale of $O(1/(1-\rho))$. Since the diffusion time scale is of the order $O(1/(1-\rho)^2)$ it is tempting to conclude that the only configurations of the customer populations that matter are configurations on the invariant manifold. These configurations depend on the workload $w$ at the CPU, which then is expected to be the driving force of randomness. The goal of the next section is to make this statement rigorous.

\section{State-space collapse in heavy traffic}
\label{sec:state-space-collapse}

We are now ready to develop a diffusion approximation for the process describing the number of customers in the system, which we sometimes also refer to as the head-count process.
Consider a sequence of such processes indexed by $n$.
As $n\to\infty$,
$
  \lambda^n\to\lambda
$.
Let $\gamma^n=(I-P^T)^{-1}\lambda^n$, and
\begin{equation*}
%  \label{eq:rho-n}
  \rho^n=(\gamma^n)^T\beta.
\end{equation*}
We assume that
\begin{equation}
  \label{eq:HT}
  \rho^n=1-\theta/n>0,\textrm{ and } K^n_i = K_in.
\end{equation}
(One way to achieve this is to set $\lambda_i^n=\lambda_i(1-\theta/n)$.)
We are interested in the limit of the diffusion scaled process
\begin{equation*}
%  \label{eq:x-diffusion}
  \ds X(t)=\frac{1}{n}X^n(n^2t)
\end{equation*}
as $n\rightarrow\infty$, in which case the system approaches heavy traffic. It turns out that the choice $K^n_i = K_in$ gives rise to a limit model in which the fraction of time the system is congested is non-trivial (i.e.\ strictly between 0 and 1). For example, in the single-node case, the results in \cite{ZhangZwart2008} imply that the time-dependent delay probability $P(X^n_1(n^2t) >K_1^n)$, as well as the stationary delay probability $P(X^n_1(\infty) >K_1^n)$, converge to a quantity between 0 and 1. This enables one to obtain non-trivial and explicit approximations of the delay probability.

A starting point of our analysis is to recall the well-known (see e.g.\ \cite{Asmussen2003}) heavy-traffic limit theorem for the workload process at the CPU layer. Let $\ds W(t)= W^{n}(n^2t)/n, t\geq 0$ be the scaled workload process. Then $\ds W(\cdot)$ converges to a reflected Brownian motion (RBM) $W^*(\cdot)$ with drift $-\theta$ and variance $\sigma^2=\E(S)(1+c_s^2)=E[S^2]/E[S]$, where $c_s^2=Var(S)/\E^2(S)$. According to the calculation in Section~\ref{sec:model-description}, $\sigma^2= {a^T \tau^{(2)}}/{(a^T\tau)}$.

Our main result is that $\ds X (t)$ converges to a process that can be described completely  in terms of $W^*(t)$, using the insights developed for the critical fluid model in the previous section.
To this end, define the map $\Delta:\R_+\to\R_+^J$ by
\begin{align}
  \label{eq:lift-map-1}
  \Delta_1(w)&=
  \frac{w\wedge w^*}{w^*}K_1
  +\frac{(w-w^*)^+}{\tau_1},\\
  \label{eq:lift-map-i}
  \Delta_i(w)&=
  \frac{w\wedge w^*}{w^*}\frac{\mu_1K_1}{\gamma_1}\frac{\gamma_i}{\mu_i},
  \quad i=2,\ldots, J.
\end{align}

This map is called \textsl{lifting map}, as it will be used to construct the multi-dimensional limiting queue length process from the one-dimensional limiting workload process. In fact, our next result is a consequence of the fact that, in heavy traffic, $\ds X(t) \approx \Delta(\ds W(t))$, and making this statement rigorous is in fact a key ingredient of the proof, which is based on state-space-collapse techniques as developed by Bramson~\cite{Bramson1998}. Again, the proof is omitted because of space limitations.

\begin{thm}[Diffusion limit]\label{thm:diffusion-limit}
Suppose that $X^n(0)=0$ for all $n$. Then, the diffusion-scaled process $\ds X$ converges weakly to the limit $X^*$ in heavy traffic. The limit $X^*$ can be characterized as $X^*(t)=\Delta(W^*(t))$, i.e.\
  \begin{align*}
%    \label{eq:SSC-1}
    X^*_1(t)&=\frac{W^*(t)\wedge w^*}{w^*}K_1+\frac{(W^*(t)-w^*)^+}{\tau_1},\\
 %   \label{eq:SSC-i}
    X^*_i(t)&=\frac{W^*(t)\wedge w^*}{w^*}\frac{\mu_1K_1}{\gamma_1}\frac{\gamma_i}{\mu_i},
  \quad i=2,\ldots, J.
  \end{align*}
\end{thm}

Note the similarity of the lifting map and the quantities $x^\dagger_i(w)$ that are used to define the invariant points of the critical fluid model. In fact, this is the key physical insight that justifies the title of this paper. Namely, fix $t$, take $\theta=1, n=1/(1-\rho)$ and recall that the workload fluctuates at  the time scale of $n^2$ as we set $\ds W(t)= W^{n}(n^2t)/n$. Hence, between time $tn^2$ and $tn^2+n$ the scaled workload hardly changes for $n$ large. Namely, it will be approximately $W^*(t)$ throughout this time. During this time, by the convergence result of the fluid limit presented in the previous section, $\frac{1}{n}X^n(n^2t+n)$ will have converged to $\Delta(W^*(t))$. Thus, in heavy traffic, fluctuations of the system at the layer of the individual servers occur at a much faster time scale than fluctuations at the CPU layer. If we wish to study fluctuations of the servers, we can keep the total workload at the CPU layer fixed, and if we wish to study performance of the system on the time-scale of the CPU layer, we can assume that jobs at the servers live on the invariant manifold; any deviations away from the invariant manifold will have averaged out.

We believe that these physical insights are interesting, and may also occur in other layered systems. In the next sections, we show how these insights lead to explicit and accurate approximations of the layered system under consideration.

\section{Steady-state performance approximations}
\label{s:Steady-state performance approximations}

In the previous section we have considered a sequence of systems approaching heavy traffic. The goal of the present
section is to utilize Theorem \ref{thm:diffusion-limit} and obtain performance approximations for the steady-state distribution for the original system.
The first step is to establish a heavy-traffic limit for the sequence of steady-state distributions indexed by $n$. It is well-known that the normalized steady-state workload of an $M/G/1$ queue in heavy traffic converges to an exponentially distributed random variable; i.e.\ if we consider the sequence of systems introduced in the previous section, let $W^n (\infty)$ be the steady-state workload in the $n$th system and $\ds W(\infty)=\frac{1}{n}W(\infty)$, then
\begin{equation*}
\ds W(\infty) \dto W^*(\infty),
\end{equation*}
where $\dto$ means convergence in distribution and $W^*(\infty)$ is an exponentially distributed random variable with mean $m=\frac{\sigma^2}{2\theta}$, by the classical steady-state analysis of RBM \cite{Asmussen2003}.

Since $W^*(\infty)$ can also be seen as the limit (in distribution) of $W^*(t)$ as $t\rightarrow\infty$, it is natural to expect that the heavy-traffic ($n\rightarrow\infty$) and steady-state limits ($t\rightarrow\infty$) can be interchanged when considering $\ds X(t)$. It is possible to do this in the same way as has been carried out in the single-node case \cite{ZhangZwart2008}; detailed are omitted due to space limits.
%, where uniform bounds (in $n$) on coupling times for (diffusion-scaled) work-conserving systems have been established, allowing a limit interchange to take place.
We can exploit this to derive a heavy-traffic limit theorem for $X^n(\infty)$, which is a $J$-dimensional random vector denoting the customer population in steady state in the $n$th system. Since $\Delta$ is continuous, we have the following result by the continuous mapping theorem:
\begin{equation*}
\ds X(\infty) \Rightarrow X^*(\infty) := \Delta(W^*(\infty)).
\end{equation*}
Note that
$%
%\begin{equation*}
%  \label{eq:prob-X}
  \pr(X^*_i(\infty)>x)=\pr\left(\Delta_i(W^*(\infty))>x\right).
%\end{equation*}
$

Since the distribution of $W^*(\infty)$ is explicit, as is the mapping $\Delta$, the above formula is explicit. Thus, we can develop explicit approximations for the original system that will be accurate in heavy traffic.

Recall we called our sequence of systems indexed by $n$ `a sequence of virtual systems'. The total load in the $n$th virtual system is $\rho^n=1-\theta/n$ and number of servers at node $i$ are $K_i^n$. In practice, one would like to get back to the original system, so we need to determine which virtual system is appropriate. If we take $\theta=1$, then we should take $n^*=1/(1-\rho)$, which also implies that in the fluid model the number of active servers at node $i$  should be equal to $(1-\rho)K_i$.

For our running example, the tandem network with 10 active servers at the first node, 20 active servers at the second node, and a total system load of 0.8, then $n=5$, and the relevant fluid model is the one where the number of active servers at the first node equals $2$ and at the second node $4$.

In what follows, the quantities $K_i$ represent the number of servers at node $i$ in the actual system.
The critical workload level $w^*$ can be rewritten as
\[
w^* = (1-\rho)\sum_j\rho_j\tau_jK_1/\rho_1.
\]
The right-hand side can be simplified further using \cite[Corollary III.5.3]{Asmussen2003}:
\[
w^* = (1-\rho)\sum_j\rho_j\tau_jK_1/\rho_1 = (1-\rho)K_1 \rho m/\rho_1.
\]
As $W^*$ is exponential with mean $m$, the heavy-traffic approximation of the delay probability at the bottleneck becomes
\begin{equation}%keep number regardless
 P(W^*>w^*) = e^{-(1-\rho)K_1 \frac{\rho}{\rho_1}}  \approx \rho^{K_1 \frac{\rho}{\rho_1}}=:p_d.
\end{equation}
In the second equation we used that $e^{-(1-\rho)}\approx \rho$ to obtain an approximation  more in line with the single-node approximation proposed by \cite{Avi-ItzhakHalfin1988}. Due to lack of space, we focus on one additional performance measure only, namely the expected total response time (i.e.\ the sojourn time) $E[V]$ of an arbitrary job which can be computed using Little's law:
\[
E[V] = E[\sum_j X_j]/\lambda^o \approx \frac{1/\lambda^o}{1-\rho} E[\sum_j \Delta_j (W^*)].
\]
Straightforward computations, combined with the above approximations, yield
\[
E[\sum_j \Delta_j (W^*)] \approx  (1-p_d) +p_d \frac{m}{\tau_1}.
\]
It makes sense to multiply the right-hand side with $\rho$ to obtain a result that is exact for the single-node case, and from a heavy-traffic point of view, $(\rho\approx 1)$ this still yields asymptotically accurate estimates. Putting everything together, our heavy-traffic approximation for $E[V]$ becomes
\begin{equation}\label{finalapproximation}
E[V] \approx \frac{E[S]}{1-\rho}  \left[(1-p_d) +p_d \frac{m}{\tau_1}\right].
\end{equation}
In the single node case for exponential job sizes, we have that $m=E[S]=\tau_1$ so our approximation indeed reduces to $E[S]/(1-\rho)$ which is the expected sojourn time in an $M/M/1$ queue. We now develop an extension valid for more general service times combining the insights of the heavy-traffic analysis of our network model with available results for the single node case.

\section{Extension to general job sizes}
\label{sec:extens-gener-job}

For the single-node case, Poisson arrivals, and general service times, \cite{Avi-ItzhakHalfin1988} proposed the approximation $p_d = \rho^{K_1}$ and
\begin{equation}
\label{halfin}
E[V] = (1-p_d) \frac{E[S]}{1-\rho}+ p_d \frac m{1-\rho},
\end{equation}
where, as before $m=E[S^2]/(2E[S])$. This approximation is exact for both FIFO ($K_1=1$) and PS ($K_1=\infty$), and \cite{ZhangZwart2008} shows the approximation is asymptotically exact in heavy traffic, using the same scaling procedure as in the present paper.
Note further that for $J=1$ we have that $\tau_1=E[S]$ and $\rho_1=\rho$ so \eqref{finalapproximation} and \eqref{halfin} coincide.

These considerations suggest that the approximation of $E[V]$ given in \eqref{finalapproximation} is still accurate for general service times assuming station 1 is the single bottleneck and keeping $\rho^{K_1 \frac{\rho}{\rho_1}}$.

Proving this necessitates an extension of the measure-valued framework in \cite{ZDZ2011}, which is beyond the scope of this paper. Instead, we validate our approximation with some simulation results for the two-node tandem case.

Let $\beta_i^e=\beta_i^{(2)}/2\beta_i$ be the mean residual service time of a job at station $i$. For the two-node tandem case we have $\gamma_1=\gamma_2=\lambda_1$.

Since we fixed the topology of the network we will no longer assume that node 1 is always the bottleneck.
Observing that $p_d \approx \rho^{\rho K_{i^*}/\rho_{i^*}}$ if node $i^*$ is the bottleneck, we obtain
\[
E[V] \approx \frac{E[S]}{1-\rho}   \left[(1-p_d) +p_d m \frac{1}{\tau_{i^*}}\right].
\]
The constant $m$ can be computed by noting that
\[
m= E[S^2]/(2E[S])= \frac{\rho_1}{\rho} (\beta_1^e+\beta_2)+ \frac{\rho_2}{\rho} \beta_2^e.
\]
We now present some numerical results for the case that both service times follow a hyper-exponential distribution.
In all examples, we focus on a moderately loaded system with $\rho=0.7$. We let the coefficient of variation of the service times range from 4 to 10 at both nodes (in fact we take the same parameters as done in the experiment of \cite{Wei04}). Note that the squared coefficient of variation $c_i^2$ of the service time at node $i$ satisfies
$c_i^2=\beta_i^{(2)}/\beta_i^2-1$.
\begin{table}
\caption{Simulation results}
\centering
\begin{tabular}{|c|c|c|}
  \hline
  $(\beta_1,\beta_2,c_1^2, c_2^2,K_1,K_2)$ & approximation & simulation \\\hline
(1, 2, 4, 4, 3, 7)      & 10.24  & 10.41  \\
(1, 2, 4, 10, 4, 6)     & 11.37  & 10.71  \\
(1, 2, 10, 4, 4,  6)    & 10.77  & 10.57  \\
(1, 2, 10, 10, 4, 6)    & 11.58  & 10.87  \\
(2, 1,  4, 4,   6, 4)   & 10.24  & 10.49  \\
(2, 1,  4, 10, 6, 4)    & 10.38  & 10.70  \\
(2, 1,  10, 4,  6, 4)   & 10.78  & 10.98  \\
(2, 1,  10, 10, 6, 4)   & 10.91  & 11.18  \\
(1, 10, 4, 4, 2, 8)     & 38.86  & 37.43  \\
(1, 10, 4, 10, 2, 8)    & 43.20  &  37.83 \\
(1, 10, 10, 4, 2, 8)    & 38.91  & 37.53  \\
(1, 10, 10, 10, 2, 8)   & 43.24  & 37.97  \\
(10, 1, 4, 4, 8, 2)     & 38.52  & 38.88  \\
(10, 1, 4, 10, 8, 2)    & 38.56  & 39.11  \\
(10, 1, 10, 4, 8, 2)    & 42.46  & 40.77  \\
(10, 1, 10, 10, 8, 2)   & 42.50  & 41.00  \\
  \hline
\end{tabular}
\end{table}

Generally, the heavy-traffic approximations are quite accurate, always within 15\% of the outcome predicted by simulation, and in several cases the error is as small as 2\%. We find that the results become less accurate if the coefficient of variation of the service time at the bottleneck is high. Similar conclusions can be drawn for higher values of the load and for larger networks.

\section{Concluding remarks}\label{s:Concluding remarks}

By establishing fluid and diffusion approximations of a two-layered queuing network, we have shown that, under critical loading, different layers in the network operate at different time scales. From the macroscopic CPU viewpoint, the system behaves as a simple one-server queue, which when critically loaded fluctuates at a time scale of $O(1/(1-\rho)^2)$. The network dynamics at the other layer evolve at a faster time scale $O(1/(1-\rho))$, thus always reaching an invariant point as if the total workload at the CPU were constant.

We have established this result by introducing fluid and diffusion approximation techniques to study layered networks. It is interesting to examine the potential of such techniques to analyze other layered networks, such as those in \cite{Rolia95,Woodside95}.

For our model, state-space collapse was established as a consequence of the single bottleneck assumption. Driven by curiosity, we are currently extending the analysis to multiple bottlenecks, although we note that the single bottleneck assumption will typically be an artefact of the fact that the buffer sizes $K_i$ need to be chosen as integers in implementations.

Another interesting topic is to allow for general job sizes, as well as time-varying arrival rates. Finally, we expect the results to be directly useful to dimension thread-pools in web servers in a static fashion. The techniques in this paper are likely to be useful for dynamic thread-pool dimensioning as well, as the application of the techniques in this paper seems promising to formulate tractable (Brownian) control problems.

\bibliographystyle{IEEEtran}
\bibliography{pub,nyp,cs}
\end{document}